\documentclass{article}

\usepackage{graphicx}
\newtheorem{thm}{Theorem}[section]

\newtheorem{de}[thm]{Definition}

\begin{document}

\title{A Direct Prediction of the Shape Parameter---a purely scattered data approach}         
\author{Lin-Tian Luh\\Department of Financial and Computational Mathematics\\ Providence University\\
Shalu Area, Taichung City, Taiwan\\ Email:ltluh@pu.edu.tw}        
\date{}

\maketitle

{\bf Abstract.} In this paper we present an approach which predicts directly without search the optimal choice of the shape parameter $c$ contained in the multiquadrics  $(-1)^{\lceil \frac{\beta}{2} \rceil}(c^{2}+\|x\|^{2})^{\frac{\beta}{2}},\ \beta>0$, and the inverse multiquadrics $(c^{2}+\|x\|^{2})^{\frac{\beta}{2}},\ \beta <0$. Unlike the simplex scheme where the data points are required to be evenly spaced, as in a recent paper of the author, here we allow them to be arbitrarily scattered in the simplex, making it much more useful. The drawback is that its theoretical ground is not so strong as in the evenly spaced data setting. However, experiments show that it works well. The experimentally optimal value of $c$ coincides with the theoretically predicted one. Since the fill distance involved is always of reasonable size, this approach is supposed to be practically useful. 
\\
{\bf AMS}: 41A05,65D05,65M15,65M70,65N15,65N50

\section{Introduction}       
In this paper the approximated functions lie in a space called $B_{\sigma}$ as in the following definition. 

\begin{de}
For any $\sigma>0$, the class of band-limited functions $f$ in $L^{2}(R^{n})$ is
$$B_{\sigma}:=\{ f\in L^{2}(R^{n}):\ \hat{f}(\xi)=0\ if\ |\xi|>\sigma \},$$
where $\hat{f}$ denotes the Fourier transform of $f$.
\end{de}
This function space looks small. In fact, it plays only an intermediate role in the interpolation. Via the $B_{\sigma}$ functions, all functions in the Sobolev space can be interpolated by the multiquadrics or inverse multiquadrics, as will be explained further in the paper.

The radial function we adopt is
\begin{eqnarray}
  h(x):=\Gamma(-\frac{\beta}{2})(c^{2}+\|x\|^{2})^{\frac{\beta}{2}},\ \beta\in R\backslash 2N_{\geq 0},\ c>0,
\end{eqnarray}
where $\|x\|$ is the Euclidean norm of $x\in R^{n}$, $\Gamma$ is the classical gamma function, and $\beta,\ c$ are constants. Note that this definition is slightly different from the one mentioned in the abstract.The definition (1) will  simplify its Fourier transform and the presentation of  our central theorem. The function $h(x)$ in (1) is conditionally positive definite (c.p.d.) of order $m=\max \{0,\lceil \frac{\beta}{2}\rceil \}$ where $\lceil \frac{\beta}{2}\rceil $ means the smallest integer greater than or equal to $\frac{\beta}{2}$. For further details we refer the reader to Madych and Nelson \cite{MN1} and Wendland \cite{We}.

For any set of data points $(x_{j},y_{j}),\ j=1,\ldots,N$, where $X=\{ x_{1},\ldots,x_{N}\}$ is a subset of $R^{n}$ and $y_{j}$ are real or complex numbers, we can always find an interpolant of the form
\begin{eqnarray}
  s(x)=p(x)+\sum_{j=1}^{N}c_{j}h(x-x_{j}),
\end{eqnarray}
where $p(x)$ is a polynomial in $P_{m-1}^{n}, m=max\{0,\lceil \frac{\beta}{2}\rceil \}$, and $c_{j}$ are coefficients to be chosen, as long as $X$ is a determining set for $P_{m-1}^{n}$. Interested readers can find these in \cite{MN1}.

Although we are interested only in scattered data, our criteria of choosing $c$ are developed from a core theorem which involves a simplex scheme with evenly spaced data points, as mentioned in the abstract. Therefore it is necessary to make a brief description of the evenly spaced scheme.

Let $T_{n}$ denote an $n$-simplex in $R^{n}$. Then $T_{1}$ is just a line segment, $T_{2}$ is a triangle, and $T_{3}$ is a tetrahedron with four vertices. The exact definition can be found in Fleming \cite{Fl}.

Let $v_{i},\ 1\leq i\leq n+1$ be the vertices of $T_{n}$. Then any point $x\in T_{n}$ can be written as a convex combination of the vertices:
$$x=\sum_{i=1}^{n+1}c_{i}v_{i},\ \sum_{i=1}^{n+1}c_{i}=1,\ c_{i}\geq 0.$$
The numbers $c_{1},\ldots ,c_{n+1}$ are called the barycentric coordinates of $x$.

For any $n$-simplex $T_{n}$, the evenly spaced points of degree $l$ are those points whose barycentric coordinates are of the form
$$(\frac{k_{1}}{l},\frac{k_{2}}{l}, \ldots, \frac{k_{n+1}}{l}),\ k_{i}\ nonnegative\ integers\ with \ \sum_{i=1}^{n+1}k_{i}=l.$$
Obviously, the number of evenly spaced points of degree $l$ is $N=\left( \begin{array}{c}
                                                         n+l \\
                                                          n
                                                       \end{array} \right) $. It's proven in Bos \cite{Bo} that such points do form a determining set for $P_{l}^{n}$.

Before entering our core theorem, some ingredients must be explained. Each function of the form (1) induces a function space ${\cal C}_{h,m},\ m=\max \{0,\lceil \frac{\beta}{2}\rceil \}$, called native space. Also, there is a seminorm $\|f\|_{h}$ for each $f\in {\cal C}_{h,m}$. These can be found in Luh \cite{Lu1,Lu2,Lu3}, Madych and Nelson \cite{MN1,MN2} and Wendland \cite{We}. The constants $\rho$ and $\Delta_{0}$, which are usually very small positive numbers for low dimensions,  in the theorem are determined by $n$ and $\beta$. We omit their complicated definitions and refer the reader to Luh \cite{Lu4}.

The following theorem is just our core theorem. We omit its complicated proof and take it directly from \cite{Lu4}.
\begin{thm}
  Let $h$ be as in (1). For any positive number $b_{0}$, let $C=\max \left\{ \frac{2}{3b_{0}},\frac{8\rho}{c}\right\}$ and $\delta_{0}=\frac{1}{3C}$. For any n-simplex $T_{n}$ of diameter $r$ satisfying $\frac{1}{3C}\leq r\leq \frac{2}{3C}$ (note that $\frac{2}{3C}\leq b_{0}$), if $f\in {\cal C}_{h,m}$,
\begin{eqnarray}
  |f(x)-s(x)|\leq 2^{\frac{n+\beta-7}{4}}\pi^{\frac{n-1}{4}}\sqrt{n\alpha_{n}}c^{\frac{\beta}{2}}\sqrt{\Delta_{0}}\sqrt{3C}\sqrt{\delta}(\lambda')^{\frac{1}{\delta}}\|f\|_{h}
\end{eqnarray}
holds for all $x\in T_{n}$ and $0<\delta\leq\delta_{0}$, where $s(x)$ is defined as in (2) with $x_{1},\ldots ,x_{N}$ the evenly spaced points of degree $l$ in $T_{n}$ satisfying $\frac{1}{3C\delta}\leq l\leq \frac{2}{3C\delta}$. The constant $\alpha_{n}$ denotes the volume of the unit ball in $R^{n}$, and $0<\lambda'<1$ is given by 
$$\lambda'=\left(\frac{2}{3}\right)^{\frac{1}{3C}}$$
which only in some cases mildly depends on the dimension n.

\end{thm}
{\bf Remark.} Note that as the degree $l$ of the evenly spaced data points increases, the number $\delta$ will decrease, making the upper bound in (3) small. Hence $\delta$ can be regarded in spirit as the well-known fill distance. It is natural to ask what will happen if one regards $\delta$ completely the same as the fill distance. If so, the requirement that the centers $x_{1},\ldots ,x_{N}$ be evenly spaced in the simplex can be dropped, making this theorem much more useful. In fact, this is just the central idea of this paper.

\section{Criteria of choosing $c$}
The number $b_{0}$ in Theorem 1.2 controls the diameter of the domain. The upper bound in (3) is greatly related to the choice of $c$. In Luh \cite{Lu4} (3) is successfully transformed into a pleasant and lucid form which shows the influence of $c$ explicitly. There are three cases: (i)$\beta>0$ and $n\geq 1$, (ii)$\beta<0$ and $n+\beta \geq 1$, or $n+\beta=-1$, and (iii)$\beta=-1$ and $n=1$.

For (i) and (ii), we have
\begin{eqnarray}
 |f(x)-s(x)|\leq d_{0}\sigma^{\frac{1+\beta+n}{4}}MN(c)\|f\|_{L^{2}(R^{n})}
\end{eqnarray}
where $d_{0}$ is a small (for low dimensions) constant independent of $c,\ \sigma$, and $f$, and $MN(c)$ is a function of $c$ defined by
\begin{eqnarray}
MN(c):= \left\{ \begin{array}{ll}
                    \sqrt{8\rho}\cdot c^{\frac{\beta-1-n}{4}}\cdot e^{\frac{c\sigma}{2}}\cdot\left(\frac{2}{3}\right) ^{\frac{c}{24\delta\rho}}   & \mbox{if $24\rho \delta \leq c<12\rho b_{0},$} \\
                    \sqrt{\frac{2}{3b_{0}}}\cdot c^{\frac{1+\beta-n}{4}}\cdot e^{\frac{c\sigma}{2}}\cdot \left( \frac{2}{3}\right) ^{\frac{b_{0}}{2\delta}} 
                            & \mbox{if $c\geq 12\rho b_{0}$.}

                  \end{array}
          \right.
\end{eqnarray}
For (iii), we have 

\begin{eqnarray}
|f(x)-s(x)|\leq d_{0}'MN(c)\|f\|_{L^{2}(R^{n})},
\end{eqnarray}
where $d_{0}'$ is only a bit different from $d_{0}$, and $MN(c)$ is defined by

\begin{eqnarray}
MN(c):=\left\{ \begin{array}{ll}
                   \sqrt{8\rho}\cdot c^{\frac{\beta-1}{2}}\cdot \left( \frac{2}{3}\right) ^{\frac{c}{24\delta\rho}}M(c)   & \mbox{if $24\rho \delta\leq c<12\rho b_{0},$} \\
                   \sqrt{\frac{2}{3b_{0}}}\cdot c^{\frac{\beta}{2}}\cdot \left( \frac{2}{3}\right) ^{\frac{b_{0}}{2\delta}}M(c)& \mbox{if $c\geq 12\rho b_{0},$} 
                \end{array}
        \right. 
\end{eqnarray}
where $M(c)$ is defined by
$$M(c):=\left\{  \begin{array}{ll}
                  \frac{1}{\sqrt{{\cal K}_{0}(1)}}   & \mbox{if $c\leq \frac{1}{\sigma},$} \\
                  \left\{ \frac{1}{{\cal K}_{0}(1)}+2\sqrt{3}\sqrt{c\sigma}e^{c\sigma}\right\}^{1/2} & \mbox{if $c>\frac{1}{\sigma},$}
                \end{array}
        \right. $$
${\cal K}_{0}$ being the modified Bessel function, for $c\in [24\rho\delta, \infty)$.

In the following text of this section the interpolation domain is a simplex in $R^{n}$ and the parameter $\delta$ is interpreted as the well-known fill distance. For the definition of fill distance, we refer the reader to Madych and Nelson \cite{MN2} and Wendland \cite{We}. Then we have the following criteria of choosing $c$.\\
\\
{\bf Case 1.} Let $f\in B_{\sigma}, \ \sigma>0$. If (i) or (ii) holds, for any given $b_{0}>0$ and $\delta<\frac{b_{0}}{2}$, the optimal choice of $c$ in the interval $[24\rho\delta,\infty)$ for the interpolation of $f$ by $s$ defined in (2) in a simplex of diameter less than or equal to $b_{0}$ is the number minimizing $MN(c)$ in (5).\\
\\
{\bf Case 2.} Let $f\in B_{\sigma},\ \sigma>0$. If (iii) holds, for any given $b_{0}>0$ and $\delta<\frac{b_{0}}{2}$, the optimal choice of $c$ in the interval $[24\rho\delta,\infty)$ for the interpolation of $f$ by $s$ defined in (2) in a simplex of diameter less than or equal to $b_{0}$ is the number minimizing $MN(c)$ in (7).\\
\\
The number $\rho$ in this paper is always equal or close to 1 and $24\rho\delta$ is usually very small. Furthermore, experiments show that the optimal value of $c$ never falls into the interval $(0,24\rho\delta)$. Hence we have essentially dealt with $c\in(0,\infty)$. The relaxation of $\delta$ from its original definition to fill distance is natural and reasonable since in Theorem 1.2 it behaves in spirit exactly the same as the fill distance. As for the shape of the domain, we do not know how important it is. Maybe more experimental evidences should be collected first. For now, it does not seem to be possible to get rid of the simplex requirement in Theorem 1.2, both in theory and practice.

\section{Experiments}
We provide two sets of experiments here. Although we concern ourselves mainly with the scattered data setting, as a comparison, the evenly spaced data setting is also tested.
\subsection{The evenly spaced data setting}
Let us investigate Case 2. of the last section, i.e., $\beta=-1$ and $n=1$. Suppose $\sigma=1,b_{0}=5$. The curves of the MN function $MN(c)$ are presented in Figures 1-5, where $\delta$ was defined in Theorem 1.2.
\begin{figure}[h]
\centering
\includegraphics[scale=1.0]{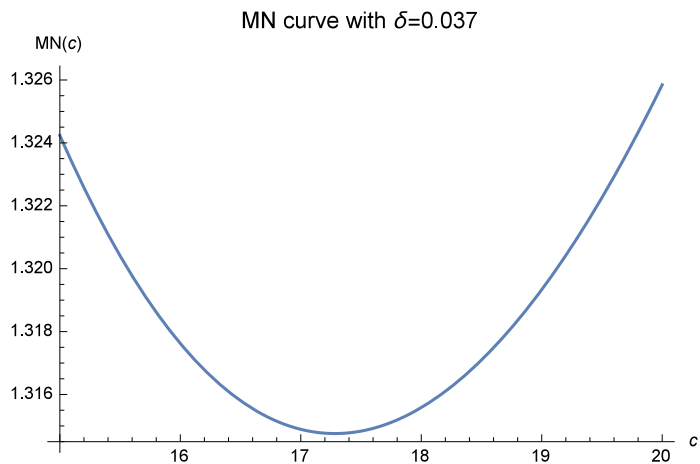}
\caption{Here $n=1,\ \beta=-1, b_{0}=5$ and $\sigma=1$.}

\end{figure}

\clearpage

\begin{figure}[t]
\centering
\includegraphics[scale=1.0]{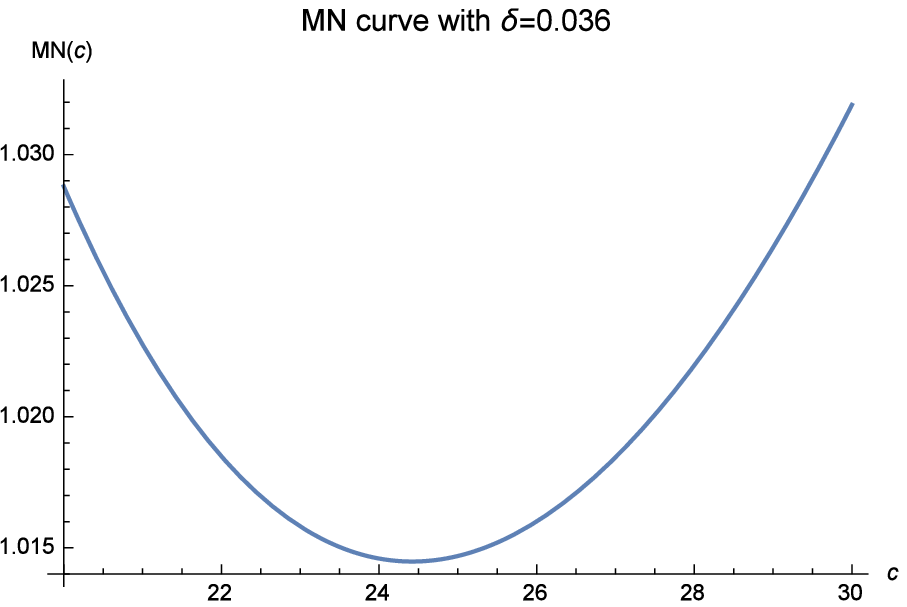}
\caption{Here $n=1,\ \beta=-1, b_{0}=5$ and $\sigma=1$.}

\vspace{2cm}

\centering
\includegraphics[scale=1.0]{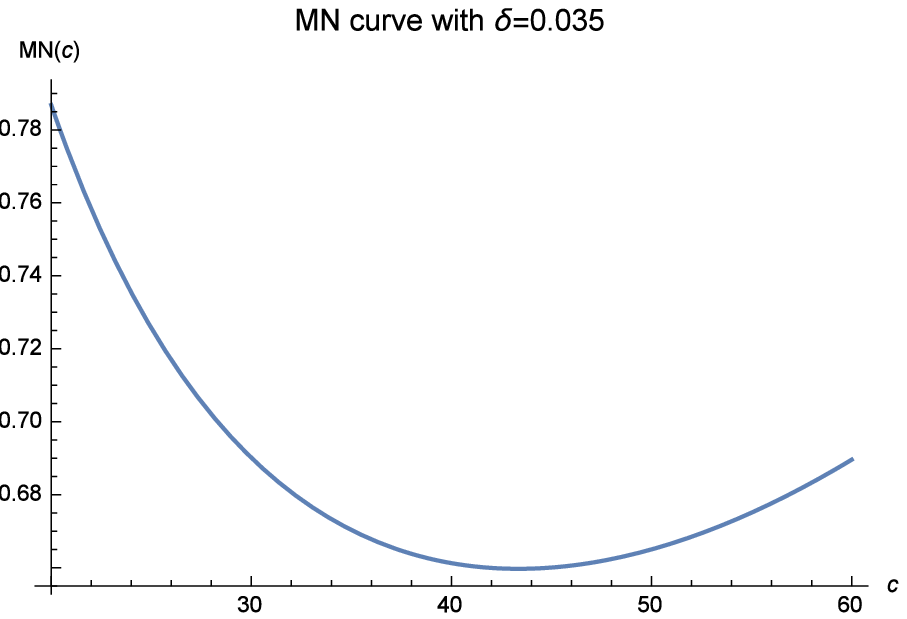}
\caption{Here $n=1,\ \beta=-1, b_{0}=5$ and $\sigma=1$.}

\end{figure}

\clearpage

\begin{figure}[t]
\centering
\includegraphics[scale=1.0]{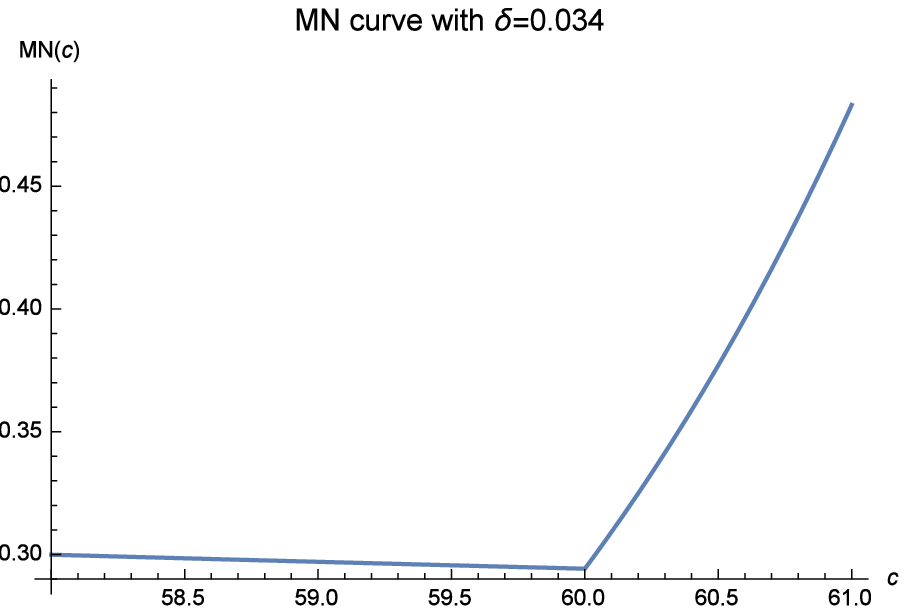}
\caption{Here $n=1,\ \beta=-1, b_{0}=5$ and $\sigma=1$.}

\vspace{2cm}

\centering
\includegraphics[scale=1.0]{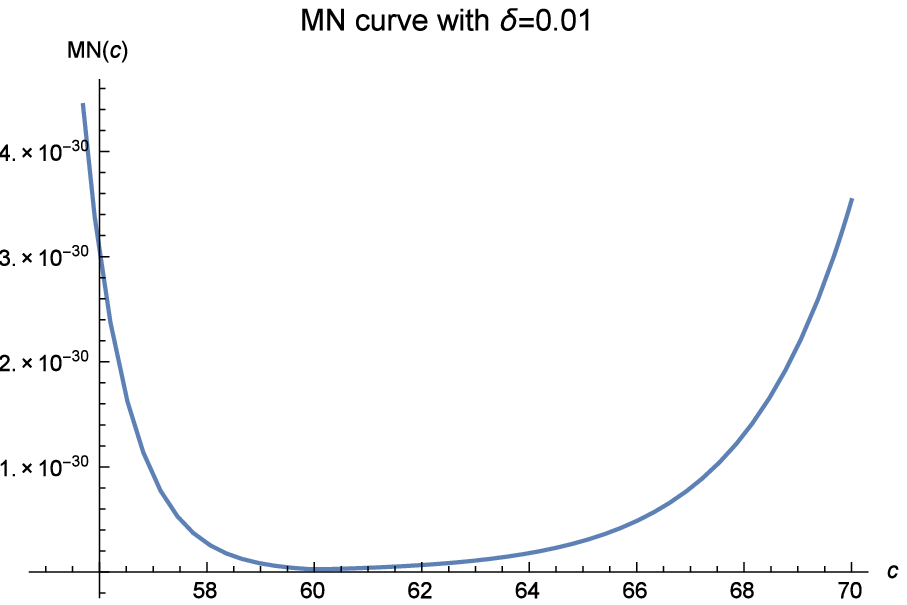}
\caption{Here $n=1,\ \beta=-1, b_{0}=5$ and $\sigma=1$.}

\end{figure}

\clearpage

In Figures 1-5, one easily finds that as $\delta$ decreases, the optimal values of $c$ move rapidly to 60. It strongly suggests that one should choose $c=60$ as the optimal value. Now we can start our experiment.

In this experiment the approximated function adopted is
$$ f(x):= \left\{ \begin{array}{ll}
                     \frac{\sin x}{x}   & \mbox{if $x\neq 0$} \\
                            \           \     1   & \mbox{if $x=0.$}
                    \end{array}
          \right. $$                        
It is easy to check that $f\in B_{\sigma}$ for $\sigma=1$. We use $s(x)$ defined in (2) to interpolate $f(x)$ in the interval $[0,5]$. However, for simplicity, the radial function used is the one mentioned in the abstract, rather than that of (1). The numbers of the centers (interpolation points) and test points are denoted by $N_{d}$ and $N_{t}$, respectively. The centers $x_{1},\ldots,x_{N_{d}}$ are evenly spaced in $[0,5]$, and so are the test points $z_{1},\ldots,z_{N_{t}}$. We use the root-mean-square error $RMS$ to evaluate the closeness of the approximation and define
$$RMS:=\sqrt{\frac{\sum_{j=1}^{N_{t}}|f(z_{j})-s(z_{j})|^{2}}{N_{t}}}.$$
The condition number of the interpolation matrix is denoted by $COND$. As is well known, the condition numbers in the RBF interpolation are usually very large. The problem of ill-conditioning is overcome by adopting enough effective digits to the right of the decimal point, with the help of the arbitrarily precise computer software Mathematica. For example, if the condition number is $10^{150}$, we adopt at least 200 effective digits for each step of the computation. Whenever keeping 250, 300, or even more effective digits, the final results are exactly the same, it means that the ill-conditioning has been completely controlled. Therefore our results should be reliable.

There is a logical problem in our approach. According to Theorem 2.1, one should choose $c$ before determining the other parameters. However, we do not know in advance the optimal choice of $c$. Hence we fix $b_{0}, \sigma,n,\beta$, and $\delta$ first. Then the optimal $c$ can be predicted by the curves of the $MN(c)$. Once $c$ is chosen, we begin to arrange the centers according to Theorem 2.1. Here we let $l=\lfloor\frac{2}{3C\delta}\rfloor$. The results are presented in Tables 1-7.

\begin{table}[t]
\caption{$\delta=0.44$}
\centering
\tiny
\begin{tabular}{c llllll}\\[2ex]
\hline\hline \\ [1ex]
\large $c$ & \large $20$ & \large $30$ & \large $40$ & \large $50$ & \large $60$ & \large $70$  \\ [1ex]
\hline\\[1ex]

\large $RMS$ & \normalsize $ 1.5\cdot 10^{-2} $ & \normalsize $6.4\cdot 10^{-4} $ & \normalsize $ 1.9\cdot 10^{-4} $ & \normalsize $ 4.0\cdot 10^{-7} $ & \normalsize $6.4\cdot 10^{-9} $  & \normalsize $7.2\cdot 10^{-9}$   \\ [1ex]
\hline\\[1ex]

\large $COND$ & \normalsize $2.1\cdot 10^{6}$ & \normalsize $4.2\cdot 10^{12}$ & \normalsize $4.3 \cdot 10^{19}$ & \normalsize $1.4 \cdot 10^{27}$ & \normalsize $1.1 \cdot 10^{35}$ & \normalsize $3.2 \cdot 10^{36}$ \\ [1ex]
\hline\\[1ex]

\large $N_{d}$ & \normalsize $4$  & \normalsize $6$  & \normalsize $8$  & \normalsize $10$ & \normalsize $12$  & \normalsize $12$\\[1ex]
\hline\\[1ex]

\large $N_{t}$ & \normalsize $50$  & \normalsize $50$  & \normalsize $50$  & \normalsize $50$ & \normalsize $50$  & \normalsize $50$  \\[1ex]

\hline\hline \\[1ex]

\large $c$ & \large $80$ & \large $90$ & \large $100$   \\ [1ex]
\hline\\[1ex]

\large $RMS$ & \normalsize $ 7.7\cdot 10^{-9} $ & \normalsize $8.1\cdot 10^{-9} $ & \normalsize $ 8.4\cdot 10^{-9} $      \\ [1ex]
\hline \\[1ex]

\large $COND$ & \normalsize $6.0\cdot 10^{37}$ & \normalsize $7.9\cdot 10^{38}$ & \normalsize $8.1 \cdot 10^{39}$  \\ [1ex]
\hline \\[1ex]

\large $N_{d}$ & \normalsize $12$   & \normalsize $12$  & \normalsize $12$ \\[1ex]
\hline \\[1ex]

\large $N_{t}$ & \normalsize $50$  & \normalsize $50$  & \normalsize $50$  \\[1ex]
\hline\\[1ex]

\end{tabular}
\label{001}
\end{table}

\begin{table}[h]
\caption{$\delta=0.4$}
\centering
\tiny
\begin{tabular}{c llllll}\\[2ex]
\hline\hline \\ [1ex]
\large $c$ & \large $20$ & \large $30$ & \large $40$ & \large $50$ & \large $60$ & \large $70$  \\ [1ex]
\hline\\[1ex]

\large $RMS$ & \normalsize $ 4.2\cdot 10^{-3} $ & \normalsize $1.3\cdot 10^{-4} $ & \normalsize $ 3.0\cdot 10^{-6} $ & \normalsize $ 5.2\cdot 10^{-8} $ & \normalsize $6.9\cdot 10^{-10} $  & \normalsize $8.1\cdot 10^{-10}$   \\ [1ex]
\hline\\[1ex]

\large $COND$ & \normalsize $4.3\cdot 10^{8}$ & \normalsize $1.9\cdot 10^{15}$ & \normalsize $3.4 \cdot 10^{22}$ & \normalsize $1.7 \cdot 10^{30}$ & \normalsize $1.9 \cdot 10^{38}$ & \normalsize $7.7\cdot 10^{39}$ \\ [1ex]
\hline\\[1ex]

\large $N_{d}$ & \normalsize $5$  & \normalsize $7$  & \normalsize $9$  & \normalsize $11$ & \normalsize $13$  & \normalsize $13$\\[1ex]
\hline\\[1ex]

\large $N_{t}$ & \normalsize $150$  & \normalsize $150$  & \normalsize $150$  & \normalsize $150$ & \normalsize $150$  & \normalsize $150$  \\[1ex]

\hline\hline \\[1ex]

\large $c$ & \large $80$ & \large $90$ & \large $100$  \\ [1ex]
\hline\\[1ex]

\large $RMS$ & \normalsize $ 9.0\cdot 10^{-10} $ & \normalsize $9.6\cdot 10^{-10} $ & \normalsize $ 1.0\cdot 10^{-9} $     \\ [1ex]
\hline \\[1ex]

\large $COND$ & \normalsize $1.9\cdot 10^{41}$ & \normalsize $3.2\cdot 10^{42}$ & \normalsize $4.0 \cdot 10^{43}$   \\ [1ex]
\hline \\[1ex]

\large $N_{d}$ & \normalsize $13$   & \normalsize $13$  & \normalsize $13$\\[1ex]
\hline \\[1ex]

\large $N_{t}$ & \normalsize $150$  & \normalsize $150$  & \normalsize $150$ \\[1ex]
\hline\\[1ex]

\end{tabular}
\label{001}
\end{table}

\clearpage

\begin{table}[t]
\caption{$\delta=0.36$}
\centering
\tiny
\begin{tabular}{c llllll}\\[2ex]
\hline\hline \\ [1ex]
\large $c$ & \large $20$ & \large $30$ & \large $40$ & \large $50$ & \large $60$ & \large $70$  \\ [1ex]
\hline\\[1ex]

\large $RMS$ & \normalsize $ 4.2\cdot 10^{-3} $ & \normalsize $1.3\cdot 10^{-4} $ & \normalsize $ 3.3\cdot 10^{-7} $ & \normalsize $ 5.3\cdot 10^{-9} $ & \normalsize $6.6\cdot 10^{-11} $  & \normalsize $8.0\cdot 10^{-11}$   \\ [1ex]
\hline\\[1ex]

\large $COND$ & \normalsize $4.3\cdot 10^{8}$ & \normalsize $1.9\cdot 10^{15}$ & \normalsize $2.5 \cdot 10^{25}$ & \normalsize $2.0\cdot 10^{33}$ & \normalsize $3.2 \cdot 10^{41}$ & \normalsize $1.7 \cdot 10^{43}$ \\ [1ex]
\hline\\[1ex]

\large $N_{d}$ & \normalsize $5$  & \normalsize $7$  & \normalsize $10$  & \normalsize $12$ & \normalsize $14$  & \normalsize $14$\\[1ex]
\hline\\[1ex]

\large $N_{t}$ & \normalsize $150$  & \normalsize $150$  & \normalsize $150$  & \normalsize $150$ & \normalsize $150$  & \normalsize $150$  \\[1ex]

\hline\hline \\[1ex]

\large $c$ & \large $80$ & \large $90$ & \large $100$   \\ [1ex]
\hline\\[1ex]

\large $RMS$ & \normalsize $ 9.0\cdot 10^{-11} $ & \normalsize $9.7\cdot 10^{-11} $ & \normalsize $ 1.0\cdot 10^{-10} $      \\ [1ex]
\hline \\[1ex]

\large $COND$ & \normalsize $5.6\cdot 10^{44}$ & \normalsize $1.2\cdot 10^{46}$ & \normalsize $1.8 \cdot 10^{47}$   \\ [1ex]
\hline \\[1ex]

\large $N_{d}$ & \normalsize $14$   & \normalsize $14$  & \normalsize $14$\\[1ex]
\hline \\[1ex]

\large $N_{t}$ & \normalsize $150$  & \normalsize $150$  & \normalsize $150$  \\[1ex]
\hline\\[1ex]

\end{tabular}
\label{001}
\end{table}

\begin{table}[h]
\caption{$\delta=0.32$}
\centering
\tiny
\begin{tabular}{c llllll}\\[2ex]
\hline\hline \\ [1ex]
\large $c$ & \large $20$ & \large $30$ & \large $40$ & \large $50$ & \large $60$ & \large $70$  \\ [1ex]
\hline\\[1ex]

\large $RMS$ & \normalsize $ 5.1\cdot 10^{-4} $ & \normalsize $1.5\cdot 10^{-5} $ & \normalsize $ 3.8\cdot 10^{-8} $ & \normalsize $ 4.6\cdot 10^{-11} $ & \normalsize $4.7\cdot 10^{-13} $  & \normalsize $6.5\cdot 10^{-13}$   \\ [1ex]
\hline\\[1ex]

\large $COND$ & \normalsize $7.6\cdot 10^{10}$ & \normalsize $7.7\cdot 10^{17}$ & \normalsize $2.0 \cdot 10^{28}$ & \normalsize $2.8 \cdot 10^{39}$ & \normalsize $9.4 \cdot 10^{47}$ & \normalsize $9.6 \cdot 10^{49}$ \\ [1ex]
\hline\\[1ex]

\large $N_{d}$ & \normalsize $6$  & \normalsize $8$  & \normalsize $11$  & \normalsize $14$ & \normalsize $16$  & \normalsize $16$\\[1ex]
\hline\\[1ex]

\large $N_{t}$ & \normalsize $150$  & \normalsize $150$  & \normalsize $150$  & \normalsize $150$ & \normalsize $150$  & \normalsize $150$  \\[1ex]

\hline\hline \\[1ex]

\large $c$ & \large $80$ & \large $90$ & \large $100$  \\ [1ex]
\hline\\[1ex]

\large $RMS$ & \normalsize $ 7.8\cdot 10^{-13} $ & \normalsize $8.9\cdot 10^{-13} $ & \normalsize $ 9.6\cdot 10^{-13} $     \\ [1ex]
\hline \\[1ex]

\large $COND$ & \normalsize $5.2\cdot 10^{51}$ & \normalsize $1.8\cdot 10^{53}$ & \normalsize $4.2 \cdot 10^{54}$  \\ [1ex]
\hline \\[1ex]

\large $N_{d}$ & \normalsize $16$   & \normalsize $16$  & \normalsize $16$\\[1ex]
\hline \\[1ex]

\large $N_{t}$ & \normalsize $150$  & \normalsize $150$  & \normalsize $150$  \\[1ex]
\hline\\[1ex]

\end{tabular}
\label{001}
\end{table}

\clearpage

\begin{table}[t]
\caption{$\delta=0.28$}
\centering
\tiny
\begin{tabular}{c llllll}\\[2ex]
\hline\hline \\ [1ex]
\large $c$ & \large $20$ & \large $30$ & \large $40$ & \large $50$ & \large $60$ & \large $70$  \\ [1ex]
\hline\\[1ex]

\large $RMS$ & \normalsize $ 5.1\cdot 10^{-4} $ & \normalsize $2.0\cdot 10^{-6} $ & \normalsize $ 3.5\cdot 10^{-9} $ & \normalsize $ 3.4\cdot 10^{-12} $ & \normalsize $2.3\cdot 10^{-15} $  & \normalsize $3.8\cdot 10^{-15}$   \\ [1ex]
\hline\\[1ex]

\large $COND$ & \normalsize $7.6\cdot 10^{10}$ & \normalsize $3.5\cdot 10^{20}$ & \normalsize $1.5 \cdot 10^{31}$ & \normalsize $3.4 \cdot 10^{42}$ & \normalsize $2.8 \cdot 10^{54}$ & \normalsize $5.3 \cdot 10^{56}$ \\ [1ex]
\hline\\[1ex]

\large $N_{d}$ & \normalsize $3$  & \normalsize $9$  & \normalsize $12$  & \normalsize $15$ & \normalsize $18$  & \normalsize $18$\\[1ex]
\hline\\[1ex]

\large $N_{t}$ & \normalsize $150$  & \normalsize $150$  & \normalsize $150$  & \normalsize $150$ & \normalsize $150$  & \normalsize $150$  \\[1ex]

\hline\hline \\[1ex]

\large $c$ & \large $80$ & \large $90$ & \large $100$   \\ [1ex]
\hline\\[1ex]

\large $RMS$ & \normalsize $ 5.1\cdot 10^{-15} $ & \normalsize $6.2\cdot 10^{-15} $ & \normalsize $ 2.2\cdot 10^{-13} $      \\ [1ex]
\hline \\[1ex]

\large $COND$ & \normalsize $4.9\cdot 10^{58}$ & \normalsize $2.7\cdot 10^{60}$ & \normalsize $6.7 \cdot 10^{62}$  \\ [1ex]
\hline \\[1ex]

\large $N_{d}$ & \normalsize $18$   & \normalsize $18$  & \normalsize $18$ \\[1ex]
\hline \\[1ex]

\large $N_{t}$ & \normalsize $150$  & \normalsize $150$  & \normalsize $150$  \\[1ex]
\hline\\[1ex]

\end{tabular}
\label{001}
\end{table}

\begin{table}[h]
\caption{$\delta=0.24$}
\centering
\tiny
\begin{tabular}{c llllll}\\[2ex]
\hline\hline \\ [1ex]
\large $c$ & \large $20$ & \large $30$ & \large $40$ & \large $50$ & \large $60$ & \large $70$  \\ [1ex]
\hline\\[1ex]

\large $RMS$ & \normalsize $ 7.7\cdot 10^{-5} $ & \normalsize $1.5\cdot 10^{-8} $ & \normalsize $ 2.0\cdot 10^{-11} $ & \normalsize $ 6.9\cdot 10^{-16} $ & \normalsize $2.6\cdot 10^{-19} $  & \normalsize $8.4\cdot 10^{-19}$   \\ [1ex]
\hline\\[1ex]

\large $COND$ & \normalsize $1.5\cdot 10^{13}$ & \normalsize $6.4\cdot 10^{25}$ & \normalsize $8.6 \cdot 10^{36}$ & \normalsize $5.7 \cdot 10^{51}$ & \normalsize $1.5 \cdot 10^{64}$ & \normalsize $7.0 \cdot 10^{66}$ \\ [1ex]
\hline\\[1ex]

\large $N_{d}$ & \normalsize $3$  & \normalsize $11$  & \normalsize $14$  & \normalsize $18$ & \normalsize $21$  & \normalsize $21$\\[1ex]
\hline\\[1ex]

\large $N_{t}$ & \normalsize $150$  & \normalsize $150$  & \normalsize $150$  & \normalsize $150$ & \normalsize $150$  & \normalsize $150$  \\[1ex]

\hline\hline \\[1ex]

\large $c$ & \large $80$ & \large $90$ & \large $100$   \\ [1ex]
\hline\\[1ex]

\large $RMS$ & \normalsize $ 1.5\cdot 10^{-18} $ & \normalsize $2.1\cdot 10^{-18} $ & \normalsize $ 2.6\cdot 10^{-18} $      \\ [1ex]
\hline \\[1ex]

\large $COND$ & \normalsize $1.4\cdot 10^{69}$ & \normalsize $1.6\cdot 10^{71}$ & \normalsize $1.1 \cdot 10^{73}$   \\ [1ex]
\hline \\[1ex]

\large $N_{d}$ & \normalsize $21$   & \normalsize $21$  & \normalsize $21$\\[1ex]
\hline \\[1ex]

\large $N_{t}$ & \normalsize $150$  & \normalsize $150$  & \normalsize $150$  \\[1ex]
\hline\\[1ex]

\end{tabular}
\label{001}
\end{table}

\clearpage

\begin{table}[h]
\caption{$\delta=0.20$}
\centering
\tiny
\begin{tabular}{c llllll}\\[2ex]
\hline\hline \\ [1ex]
\large $c$ & \large $20$ & \large $30$ & \large $40$ & \large $50$ & \large $60$ & \large $70$  \\ [1ex]
\hline\\[1ex]

\large $RMS$ & \normalsize $ 2.0\cdot 10^{-7} $ & \normalsize $9.3\cdot 10^{-12} $ & \normalsize $ 1.5\cdot 10^{-15} $ & \normalsize $ 5.2\cdot 10^{-20} $ & \normalsize $7.0\cdot 10^{-25} $  & \normalsize $2.3\cdot 10^{-24}$   \\ [1ex]
\hline\\[1ex]

\large $COND$ & \normalsize $5.6\cdot 10^{17}$ & \normalsize $1.2\cdot 10^{31}$ & \normalsize $4.0 \cdot 10^{45}$ & \normalsize $1.0 \cdot 10^{61}$ & \normalsize $1.3 \cdot 10^{77}$ & \normalsize $2.1 \cdot 10^{80}$ \\ [1ex]
\hline\\[1ex]

\large $N_{d}$ & \normalsize $9$  & \normalsize $13$  & \normalsize $17$  & \normalsize $21$ & \normalsize $25$  & \normalsize $25$\\[1ex]
\hline\\[1ex]

\large $N_{t}$ & \normalsize $150$  & \normalsize $150$  & \normalsize $150$  & \normalsize $150$ & \normalsize $150$  & \normalsize $150$  \\[1ex]

\hline\hline \\[1ex]

\large $c$ & \large $80$ & \large $90$ & \large $100$  \\ [1ex]
\hline\\[1ex]

\large $RMS$ & \normalsize $ 1.0\cdot 10^{-23} $ & \normalsize $2.2\cdot 10^{-23} $ & \normalsize $ 3.3\cdot 10^{-23} $     \\ [1ex]
\hline \\[1ex]

\large $COND$ & \normalsize $1.3\cdot 10^{83}$ & \normalsize $3.6\cdot 10^{85}$ & \normalsize $5.7 \cdot 10^{87}$  \\ [1ex]
\hline \\[1ex]

\large $N_{d}$ & \normalsize $25$   & \normalsize $25$  & \normalsize $25$ \\[1ex]
\hline \\[1ex]

\large $N_{t}$ & \normalsize $150$  & \normalsize $150$  & \normalsize $150$  \\[1ex]
\hline\\[1ex]

\end{tabular}
\label{001}
\end{table}

In Tables 1-7 it is easily seen that the optimal values of $c$ are always 60, as predicted by the MN curves. Hence our approach of finding the optimal $c$ is extremely reliable for the evenly spaced data setting. What is noteworthy is that in these tables, the numbers of data points used are not always the same, for the same $\delta$. This results from the requirement of Theorem 1.2. According to Theorem 1.2, one should choose $c$ first and then arrange the centers by the value of $c$. 

\subsection{The scattered data setting}
Now we begin to test our theoretical prediction of the optimal value of $c$ when the data points are purely scattered. We use the Mathematica command Random[$\cdot$] to generate a random number between 0 and 1. The interpolation domain is still [0, 5]. The interval [0, 5] is divided into subintervals of width $\delta$. Each subinterval contains a random number. For example, if $[a_{i},b_{i}]$ is a subinterval, then $x_{i}=a_{i}+Random[\cdot]*\delta$ is an interpolation center in this subinterval. If $5/\delta$ is not an integer, then 5 is set to be the interpolation center of the rightmost subinterval. Obviously the fill distance in this setting is $\delta$. Then we use $s(x)$ in (2) with the gamma function replaced by 1 to interpolate $f(x)$ defined in subsection 3.1. The results are presented in Tables 8-14.

\clearpage

 \begin{table}[t]
\caption{$\delta=0.48$}
\centering
\tiny
\begin{tabular}{c llllll}\\[2ex]
\hline\hline \\ [1ex]
\large $c$ & \large $10$ & \large $16$ & \large $18$ & \large $20$ & \large $30$ & \large $40$  \\ [1ex]
\hline\\[1ex]

\large $RMS$ & \normalsize $ 1.4\cdot 10^{-7} $ & \normalsize $1.4\cdot 10^{-8} $ & \normalsize $ 1.2\cdot 10^{-8} $ & \normalsize $ 2.1\cdot 10^{-8} $ & \normalsize $3.8\cdot 10^{-8} $  & \normalsize $1.0\cdot 10^{-7}$   \\ [1ex]
\hline\\[1ex]

\large $COND$ & \normalsize $1.2\cdot 10^{17}$ & \normalsize $1.1\cdot 10^{21}$ & \normalsize $1.1 \cdot 10^{22}$ & \normalsize $8.6 \cdot 10^{22}$ & \normalsize $2.7 \cdot 10^{26}$ & \normalsize $8.2 \cdot 10^{28}$ \\ [1ex]
\hline\\[1ex]

\large $N_{d}$ & \normalsize $11$  & \normalsize $11$  & \normalsize $11$  & \normalsize $11$ & \normalsize $11$  & \normalsize $11$\\[1ex]
\hline\\[1ex]

\large $N_{t}$ & \normalsize $40$  & \normalsize $40$  & \normalsize $40$  & \normalsize $40$ & \normalsize $40$  & \normalsize $40$  \\[1ex]

\hline\hline \\[1ex]

\large $c$ & \large $50$ & \large $60$ & \large $70$   \\ [1ex]
\hline\\[1ex]

\large $RMS$ & \normalsize $ 1.4\cdot 10^{-7} $ & \normalsize $6.8\cdot 10^{-8} $ & \normalsize $ 7.5\cdot 10^{-8} $      \\ [1ex]
\hline \\[1ex]

\large $COND$ & \normalsize $7.1\cdot 10^{30}$ & \normalsize $1.5\cdot 10^{32}$ & \normalsize $3.3 \cdot 10^{33}$   \\ [1ex]
\hline \\[1ex]

\large $N_{d}$ & \normalsize $11$   & \normalsize $11$  & \normalsize $11$ \\[1ex]
\hline \\[1ex]

\large $N_{t}$ & \normalsize $40$  & \normalsize $40$  & \normalsize $40$  \\[1ex]
\hline\\[1ex]

\end{tabular}
\label{001}
\end{table}

\begin{table}[h]
\caption{$\delta=0.40$}
\centering
\tiny
\begin{tabular}{c llllll}\\[2ex]
\hline\hline \\ [1ex]
\large $c$ & \large $20$ & \large $25$ & \large $28$ & \large $30$ & \large $40$ & \large $50$  \\ [1ex]
\hline\\[1ex]

\large $RMS$ & \normalsize $ 1.6\cdot 10^{-10} $ & \normalsize $1.4\cdot 10^{-10} $ & \normalsize $ 1.1\cdot 10^{-10} $ & \normalsize $ 3.4\cdot 10^{-11} $ & \normalsize $5.9\cdot 10^{-10} $  & \normalsize $1.2\cdot 10^{-9}$   \\ [1ex]
\hline\\[1ex]

\large $COND$ & \normalsize $7.8\cdot 10^{27}$ & \normalsize $1.6\cdot 10^{30}$ & \normalsize $2.4 \cdot 10^{31}$ & \normalsize $1.2 \cdot 10^{32}$ & \normalsize $1.2\cdot 10^{35}$ & \normalsize $2.5\cdot 10^{37}$ \\ [1ex]
\hline\\[1ex]

\large $N_{d}$ & \normalsize $13$  & \normalsize $13$  & \normalsize $13$  & \normalsize $13$ & \normalsize $13$  & \normalsize $13$\\[1ex]
\hline\\[1ex]

\large $N_{t}$ & \normalsize $80$  & \normalsize $80$  & \normalsize $80$  & \normalsize $80$ & \normalsize $80$  & \normalsize $80$  \\[1ex]

\hline\hline \\[1ex]

\large $c$ & \large $60$ & \large $70$ & \large $90$ & \large $120$ \\ [1ex]
\hline\\[1ex]

\large $RMS$ & \normalsize $ 1.6\cdot 10^{-9} $ & \normalsize $1.8\cdot 10^{-9} $ & \normalsize $ 2.2\cdot 10^{-9} $ & \normalsize $ 2.4\cdot 10^{-9} $      \\ [1ex]
\hline \\[1ex]

\large $COND$ & \normalsize $2.0\cdot 10^{39}$ & \normalsize $7.9\cdot 10^{40}$ & \normalsize $3.3 \cdot 10^{43}$ & \normalsize $3.3 \cdot 10^{46}$   \\ [1ex]
\hline \\[1ex]

\large $N_{d}$ & \normalsize $13$   & \normalsize $13$  & \normalsize $13$ & \normalsize $13$  \\[1ex]
\hline \\[1ex]

\large $N_{t}$ & \normalsize $80$  & \normalsize $80$  & \normalsize $80$ & \normalsize $80$   \\[1ex]
\hline\\[1ex]

\end{tabular}
\label{001}
\end{table}

\clearpage

\begin{table}[t]
\caption{$\delta=0.32$}
\centering
\tiny
\begin{tabular}{c llllll}\\[2ex]
\hline\hline \\ [1ex]
\large $c$ & \large $20$ & \large $25$ & \large $30$ & \large $40$ & \large $50$ & \large $60$  \\ [1ex]
\hline\\[1ex]

\large $RMS$ & \normalsize $ 4.4\cdot 10^{-13} $ & \normalsize $2.8\cdot 10^{-13} $ & \normalsize $ 1.0\cdot 10^{-13} $ & \normalsize $ 1.3\cdot 10^{-13} $ & \normalsize $1.1\cdot 10^{-12} $  & \normalsize $2..0\cdot 10^{-12}$   \\ [1ex]
\hline\\[1ex]

\large $COND$ & \normalsize $5.8\cdot 10^{34}$ & \normalsize $4.4\cdot 10^{37}$ & \normalsize $1.0 \cdot 10^{40}$ & \normalsize $5.5\cdot 10^{43}$ & \normalsize $4.4 \cdot 10^{46}$ & \normalsize $1.0 \cdot 10^{49}$ \\ [1ex]
\hline\\[1ex]

\large $N_{d}$ & \normalsize $16$  & \normalsize $16$  & \normalsize $16$  & \normalsize $16$ & \normalsize $16$  & \normalsize $16$\\[1ex]
\hline\\[1ex]

\large $N_{t}$ & \normalsize $80$  & \normalsize $80$  & \normalsize $80$  & \normalsize $80$ & \normalsize $80$  & \normalsize $80$  \\[1ex]

\hline\hline \\[1ex]

\large $c$ & \large $70$ & \large $90$ & \large $120$   \\ [1ex]
\hline\\[1ex]

\large $RMS$ & \normalsize $ 2.7\cdot 10^{-12} $ & \normalsize $3.7\cdot 10^{-12} $ & \normalsize $ 4.5\cdot 10^{-12} $    \\ [1ex]
\hline \\[1ex]

\large $COND$ & \normalsize $1.1\cdot 10^{51}$ & \normalsize $2.0\cdot 10^{54}$ & \normalsize $1.1 \cdot 10^{58}$  \\ [1ex]
\hline \\[1ex]

\large $N_{d}$ & \normalsize $16$   & \normalsize $16$  & \normalsize $16$ \\[1ex]
\hline \\[1ex]

\large $N_{t}$ & \normalsize $80$  & \normalsize $80$  & \normalsize $80$  \\[1ex]
\hline\\[1ex]

\end{tabular}
\label{001}
\end{table}

\begin{table}[h]
\caption{$\delta=0.24$}
\centering
\tiny
\begin{tabular}{c llllll}\\[2ex]
\hline\hline \\ [1ex]
\large $c$ & \large $20$ & \large $30$ & \large $40$ & \large $48$ & \large $50$ & \large $52$  \\ [1ex]
\hline\\[1ex]

\large $RMS$ & \normalsize $ 7.5\cdot 10^{-18} $ & \normalsize $2.7\cdot 10^{-19} $ & \normalsize $ 7.3\cdot 10^{-20} $ & \normalsize $ 1.2\cdot 10^{-19} $ & \normalsize $1.0\cdot 10^{-19} $  & \normalsize $4.8\cdot 10^{-20}$   \\ [1ex]
\hline\\[1ex]

\large $COND$ & \normalsize $9.7\cdot 10^{45}$ & \normalsize $9.6\cdot 10^{52}$ & \normalsize $9.2 \cdot 10^{57}$ & \normalsize $1.3 \cdot 10^{61}$ & \normalsize $6.8\cdot 10^{61}$ & \normalsize $3.2\cdot 10^{62}$ \\ [1ex]
\hline\\[1ex]

\large $N_{d}$ & \normalsize $21$  & \normalsize $21$  & \normalsize $21$  & \normalsize $21$ & \normalsize $21$  & \normalsize $21$\\[1ex]
\hline\\[1ex]

\large $N_{t}$ & \normalsize $80$  & \normalsize $80$  & \normalsize $80$  & \normalsize $80$ & \normalsize $80$  & \normalsize $80$  \\[1ex]

\hline\hline \\[1ex]

\large $c$ & \large $54$ & \large $56$ & \large $60$ & \large $70$   \\ [1ex]
\hline\\[1ex]

\large $RMS$ & \normalsize $ 4.2\cdot 10^{-20} $ & \normalsize $1.6\cdot 10^{-19} $ & \normalsize $ 4.9\cdot 10^{-19} $ & \normalsize $ 1.6\cdot 10^{-18} $     \\ [1ex]
\hline \\[1ex]

\large $COND$ & \normalsize $1.5\cdot 10^{63}$ & \normalsize $6.3\cdot 10^{63}$ & \normalsize $9.9 \cdot 10^{64}$ & \normalsize $4.7 \cdot 10^{67}$  \\ [1ex]
\hline \\[1ex]

\large $N_{d}$ & \normalsize $21$   & \normalsize $21$  & \normalsize $21$ & \normalsize $21$  \\[1ex]
\hline \\[1ex]

\large $N_{t}$ & \normalsize $80$  & \normalsize $80$  & \normalsize $80$ & \normalsize $80$   \\[1ex]
\hline\\[1ex]

\end{tabular}
\label{001}
\end{table}

\clearpage

\begin{table}[t]
\caption{$\delta=0.20$}
\centering
\tiny
\begin{tabular}{c llllll}\\[2ex]
\hline\hline \\ [1ex]
\large $c$ & \large $20$ & \large $30$ & \large $40$ & \large $48$ & \large $50$ & \large $52$  \\ [1ex]
\hline\\[1ex]

\large $RMS$ & \normalsize $ 5.0\cdot 10^{-21} $ & \normalsize $5.1\cdot 10^{-24} $ & \normalsize $ 8.8\cdot 10^{-24} $ & \normalsize $ 5.0\cdot 10^{-24} $ & \normalsize $4.1\cdot 10^{-24} $  & \normalsize $1.5\cdot 10^{-24}$   \\ [1ex]
\hline\\[1ex]

\large $COND$ & \normalsize $3.3\cdot 10^{55}$ & \normalsize $8.1\cdot 10^{63}$ & \normalsize $7.7 \cdot 10^{69}$ & \normalsize $4.8 \cdot 10^{73}$ & \normalsize $3.4 \cdot 10^{74}$ & \normalsize $2.2 \cdot 10^{75}$ \\ [1ex]
\hline\\[1ex]

\large $N_{d}$ & \normalsize $25$  & \normalsize $25$  & \normalsize $25$  & \normalsize $25$ & \normalsize $25$  & \normalsize $25$\\[1ex]
\hline\\[1ex]

\large $N_{t}$ & \normalsize $80$  & \normalsize $80$  & \normalsize $80$  & \normalsize $80$ & \normalsize $80$  & \normalsize $80$  \\[1ex]

\hline\hline \\[1ex]

\large $c$ & \large $54$ & \large $56$ & \large $60$ & \large $70$  \\ [1ex]
\hline\\[1ex]

\large $RMS$ & \normalsize $ 1.8\cdot 10^{-24} $ & \normalsize $4.7\cdot 10^{-24} $ & \normalsize $ 7.1\cdot 10^{-24} $ & \normalsize $ 2.4\cdot 10^{-23} $      \\ [1ex]
\hline \\[1ex]

\large $COND$ & \normalsize $1.4\cdot 10^{76}$ & \normalsize $7.7\cdot 10^{76}$ & \normalsize $2.1 \cdot 10^{78}$ & \normalsize $3.4 \cdot 10^{81}$   \\ [1ex]
\hline \\[1ex]

\large $N_{d}$ & \normalsize $25$   & \normalsize $25$  & \normalsize $25$ & \normalsize $25$  \\[1ex]
\hline \\[1ex]

\large $N_{t}$ & \normalsize $80$  & \normalsize $80$  & \normalsize $80$ & \normalsize $80$   \\[1ex]
\hline\\[1ex]

\end{tabular}
\label{001}
\end{table}

\begin{table}[h]
\caption{$\delta=0.17$}
\centering
\tiny
\begin{tabular}{c llllll}\\[2ex]
\hline\hline \\ [1ex]
\large $c$ & \large $20$ & \large $30$ & \large $40$ & \large $50$ & \large $54$ & \large $56$  \\ [1ex]
\hline\\[1ex]

\large $RMS$ & \normalsize $ 1.0\cdot 10^{-26} $ & \normalsize $1.4\cdot 10^{-29} $ & \normalsize $ 1.2\cdot 10^{-30} $ & \normalsize $ 1.6\cdot 10^{-31} $ & \normalsize $2.0\cdot 10^{-31} $  & \normalsize $6.5\cdot 10^{-32}$   \\ [1ex]
\hline\\[1ex]

\large $COND$ & \normalsize $1.2\cdot 10^{66}$ & \normalsize $1.7\cdot 10^{76}$ & \normalsize $2.8 \cdot 10^{83}$ & \normalsize $1.1 \cdot 10^{89}$ & \normalsize $9.9 \cdot 10^{90}$ & \normalsize $8.1 \cdot 10^{91}$ \\ [1ex]
\hline\\[1ex]

\large $N_{d}$ & \normalsize $30$  & \normalsize $30$  & \normalsize $30$  & \normalsize $30$ & \normalsize $30$  & \normalsize $30$\\[1ex]
\hline\\[1ex]

\large $N_{t}$ & \normalsize $80$  & \normalsize $80$  & \normalsize $80$  & \normalsize $80$ & \normalsize $80$  & \normalsize $80$  \\[1ex]

\hline\hline \\[1ex]

\large $c$ & \large $58$ & \large $60$ & \large $70$ & \large $100$ & \large $120$  \\ [1ex]
\hline\\[1ex]

\large $RMS$ & \normalsize $ 7.9\cdot 10^{-32} $ & \normalsize $1.8\cdot 10^{-31} $ & \normalsize $ 2.4\cdot 10^{-31} $ & \normalsize $ 1.1\cdot 10^{-29} $ & \normalsize $3.2\cdot 10^{-29} $     \\ [1ex]
\hline \\[1ex]

\large $COND$ & \normalsize $6.2\cdot 10^{92}$ & \normalsize $4.4\cdot 10^{93}$ & \normalsize $3.3 \cdot 10^{97}$ & \normalsize $3.2\cdot 10^{106}$ & \normalsize $1.2 \cdot 10^{111}$  \\ [1ex]
\hline \\[1ex]

\large $N_{d}$ & \normalsize $30$   & \normalsize $30$  & \normalsize $30$ & \normalsize $30$  & \normalsize $30$\\[1ex]
\hline \\[1ex]

\large $N_{t}$ & \normalsize $80$  & \normalsize $80$  & \normalsize $80$ & \normalsize $80$  & \normalsize $80$ \\[1ex]
\hline\\[1ex]

\end{tabular}
\label{001}
\end{table}

\clearpage

\begin{table}[t]
\caption{$\delta=0.165$}
\centering
\tiny
\begin{tabular}{c llllll}\\[2ex]
\hline\hline \\ [1ex]
\large $c$ & \large $40$ & \large $50$ & \large $55$ & \large $60$ & \large $65$ & \large $70$  \\ [1ex]
\hline\\[1ex]

\large $RMS$ & \normalsize $ 3.4\cdot 10^{-31} $ & \normalsize $4.2\cdot 10^{-32} $ & \normalsize $ 3.7\cdot 10^{-32} $ & \normalsize $ 5.9\cdot 10^{-34} $ & \normalsize $3.2\cdot 10^{-32} $  & \normalsize $1.8\cdot 10^{-33}$   \\ [1ex]
\hline\\[1ex]

\large $COND$ & \normalsize $2.0\cdot 10^{87}$ & \normalsize $1.3\cdot 10^{93}$ & \normalsize $3.9 \cdot 10^{95}$ & \normalsize $7.1 \cdot 10^{97}$ & \normalsize $8.7 \cdot 10^{99}$ & \normalsize $7.4 \cdot 10^{101}$ \\ [1ex]
\hline\\[1ex]

\large $N_{d}$ & \normalsize $31$  & \normalsize $31$  & \normalsize $31$  & \normalsize $31$ & \normalsize $31$  & \normalsize $31$\\[1ex]
\hline\\[1ex]

\large $N_{t}$ & \normalsize $160$  & \normalsize $160$  & \normalsize $160$  & \normalsize $160$ & \normalsize $160$  & \normalsize $160$  \\[1ex]

\hline\hline \\[1ex]

\large $c$ & \large $80$ & \large $100$ & \large $120$   \\ [1ex]
\hline\\[1ex]

\large $RMS$ & \normalsize $ 5.0\cdot 10^{-32} $ & \normalsize $1.6\cdot 10^{-30} $ & \normalsize $ 5.6\cdot 10^{-30} $      \\ [1ex]
\hline \\[1ex]

\large $COND$ & \normalsize $2.2\cdot 10^{105}$ & \normalsize $1.4\cdot 10^{111}$ & \normalsize $8.0 \cdot 10^{115}$  \\ [1ex]
\hline \\[1ex]

\large $N_{d}$ & \normalsize $31$   & \normalsize $31$  & \normalsize $31$ \\[1ex]
\hline \\[1ex]

\large $N_{t}$ & \normalsize $160$  & \normalsize $160$  & \normalsize $160$  \\[1ex]
\hline\\[1ex]

\end{tabular}
\label{001}
\end{table}

As predicted by Figures 1-5, the optimal values of $c$ will move to 60 when $\delta$ decreases. This is supported by the results in Tables 8-14, where $\delta$ is interpreted as the fill distance.

\section{The failures of the MN curve approach}
As is well known, Newton's method of root-finding may fail whenever there are horizontal or nearly horizontal tangent lines. Similarly, our approach may also fail whenever there are nearly horizontal zones on the MN curves. Let us see a few MN curves first. If we further decrease the parameter $\delta$ in Figures 1-5 of Section 3, nearly horizontal zones will appear,  near the bottoms of the MN curves, as shown in Figures 6-8. Also, note that, for the same $\delta$, the root-mean-square errors in the tables of the preceding section are much smaller than the error bounds (4) and (6) essentially reflected by the MN function values, shown in Figures 1-5. It means that the error bounds are not very sharp. Once the curve is nearly horizontal at the bottom, the  optimal value of $c$ predicted by the MN curve may not be reliable. The actual optimal value may fall into the nearly horizontal zone. The longer this zone is, the less reliable the MN curve approach is. Experiments also show this, both in the evenly spaced and scattered data setings. Our experimental results are presented in Tables 15-18 and 19-22 for the two settings, respectively, where the $\delta$'s are smaller than those of Tables 7 and 14.

\begin{figure}[t]
\centering
\includegraphics[scale=1.0]{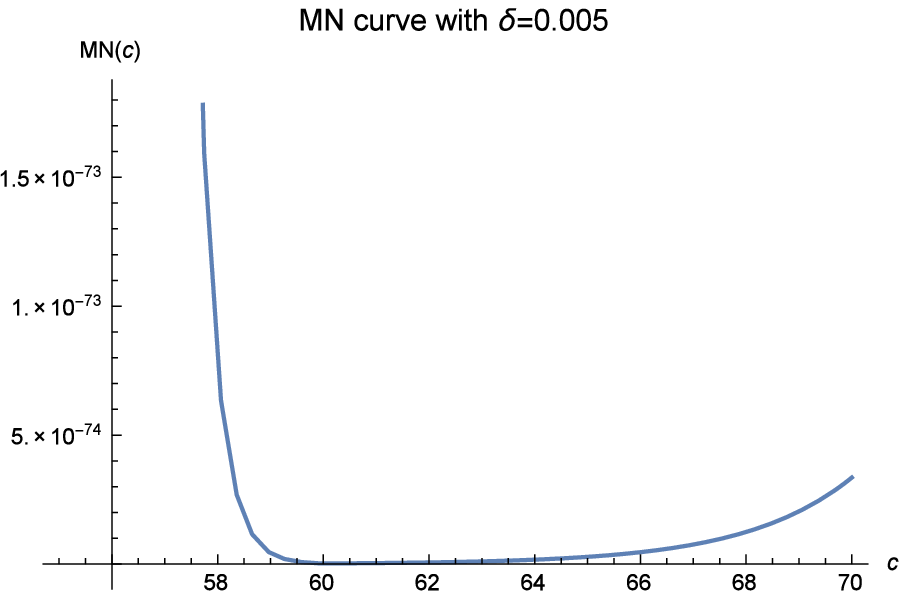}
\caption{Here $n=1,\ \beta=-1, b_{0}=5$ and $\sigma=1$.}

\vspace{2cm}

\centering
\includegraphics[scale=1.0]{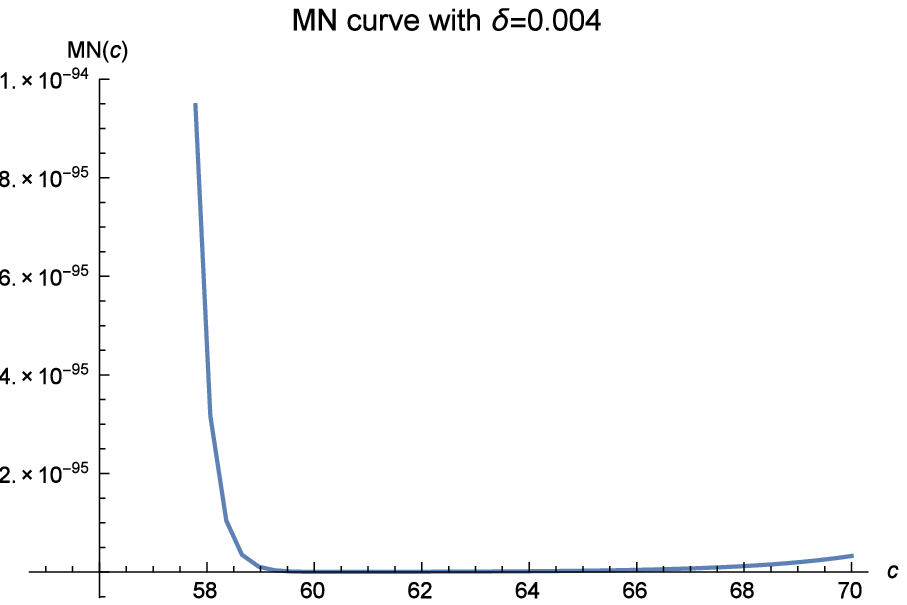}
\caption{Here $n=1,\ \beta=-1, b_{0}=5$ and $\sigma=1$.}

\end{figure}

\clearpage

\begin{figure}[t]
\centering
\includegraphics[scale=1.0]{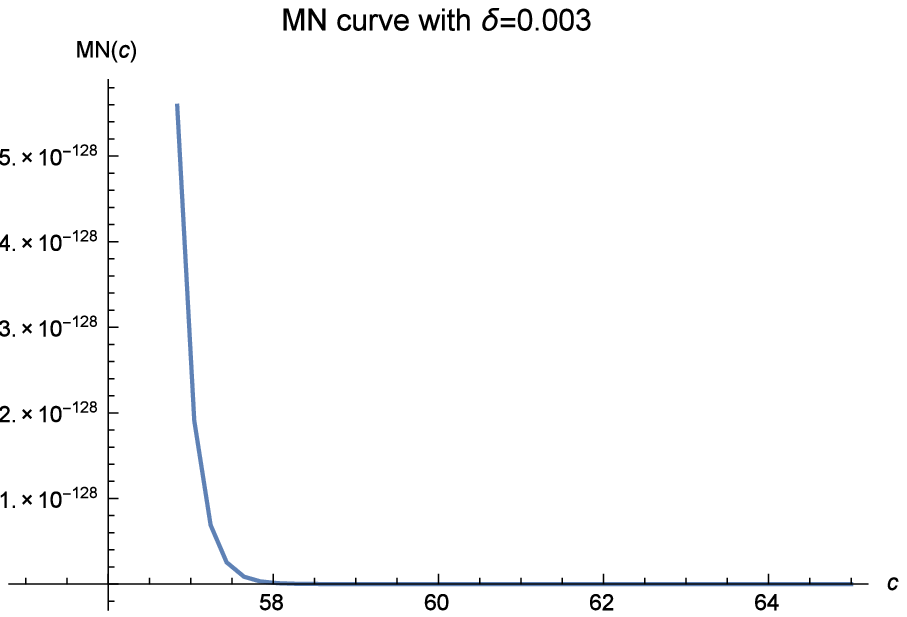}
\caption{Here $n=1,\ \beta=-1, b_{0}=5$ and $\sigma=1$.}

\end{figure}

\begin{table}[h]
\caption{$\delta=0.16$}
\centering
\tiny
\begin{tabular}{c llllll}\\[2ex]
\hline\hline \\ [1ex]
\large $c$ & \large $40$ & \large $50$ & \large $60$ & \large $70$ & \large $80$ & \large $90$  \\ [1ex]
\hline\\[1ex]

\large $RMS$ & \normalsize $ 3.8\cdot 10^{-20} $ & \normalsize $3.1\cdot 10^{-28} $ & \normalsize $ 1.9\cdot 10^{-34} $ & \normalsize $ 1.7\cdot 10^{-34} $ & \normalsize $5.4\cdot 10^{-34} $  & \normalsize $1.8\cdot 10^{-33}$   \\ [1ex]
\hline\\[1ex]

\large $COND$ & \normalsize $1.4\cdot 10^{57}$ & \normalsize $3.0\cdot 10^{79}$ & \normalsize $5.8 \cdot 10^{99}$ & \normalsize $8.2 \cdot 10^{103}$ & \normalsize $3.2 \cdot 10^{107}$ & \normalsize $4.7 \cdot 10^{110}$ \\ [1ex]
\hline\\[1ex]

\large $N_{d}$ & \normalsize $21$  & \normalsize $27$  & \normalsize $32$  & \normalsize $32$ & \normalsize $32$  & \normalsize $32$\\[1ex]
\hline\\[1ex]

\large $N_{t}$ & \normalsize $119$  & \normalsize $119$  & \normalsize $119$  & \normalsize $119$ & \normalsize $119$  & \normalsize $119$  \\[1ex]

\hline\hline \\[1ex]

\large $c$ & \large $100$ & \large $110$ & \large $120$   \\ [1ex]
\hline\\[1ex]

\large $RMS$ & \normalsize $ 1.0\cdot 10^{-32} $ & \normalsize $2.4\cdot 10^{-32} $ & \normalsize $ 4.2\cdot 10^{-32} $     \\ [1ex]
\hline \\[1ex]

\large $COND$ & \normalsize $3.2\cdot 10^{113}$ & \normalsize $1.2\cdot 10^{116}$ & \normalsize $2.6 \cdot 10^{118}$   \\ [1ex]
\hline \\[1ex]

\large $N_{d}$ & \normalsize $32$   & \normalsize $32$  & \normalsize $32$ \\[1ex]
\hline \\[1ex]

\large $N_{t}$ & \normalsize $119$  & \normalsize $119$  & \normalsize $119$  \\[1ex]
\hline\\[1ex]

\end{tabular}
\label{001}
\end{table}

\clearpage

\begin{table}[t]
\caption{$\delta=0.12$}
\centering
\tiny
\begin{tabular}{c llllll}\\[2ex]
\hline\hline \\ [1ex]
\large $c$ & \large $30$ & \large $40$ & \large $50$ & \large $60$ & \large $70$ & \large $80$  \\ [1ex]
\hline\\[1ex]

\large $RMS$ & \normalsize $ 1.4\cdot 10^{-19} $ & \normalsize $2.8\cdot 10^{-28} $ & \normalsize $ 7.4\cdot 10^{-38} $ & \normalsize $ 8.1\cdot 10^{-48} $ & \normalsize $1.1\cdot 10^{-48} $  & \normalsize $4.8\cdot 10^{-50}$   \\ [1ex]
\hline\\[1ex]

\large $COND$ & \normalsize $1.4\cdot 10^{52}$ & \normalsize $2.1\cdot 10^{77}$ & \normalsize $1.3 \cdot 10^{104}$ & \normalsize $1.4 \cdot 10^{132}$ & \normalsize $4.2 \cdot 10^{137}$ & \normalsize $2.4\cdot 10^{142}$ \\ [1ex]
\hline\\[1ex]

\large $N_{d}$ & \normalsize $21$  & \normalsize $28$  & \normalsize $35$  & \normalsize $42$ & \normalsize $42$  & \normalsize $42$\\[1ex]
\hline\\[1ex]

\large $N_{t}$ & \normalsize $150$  & \normalsize $150$  & \normalsize $150$  & \normalsize $150$ & \normalsize $150$  & \normalsize $150$  \\[1ex]

\hline\hline \\[1ex]

\large $c$ & \large $90$ & \large $100$ & \large $110$   \\ [1ex]
\hline\\[1ex]

\large $RMS$ & \normalsize $ 2.0\cdot 10^{-49} $ & \normalsize $2.7\cdot 10^{-49} $ & \normalsize $ 4.5\cdot 10^{-48} $      \\ [1ex]
\hline \\[1ex]

\large $COND$ & \normalsize $3.7\cdot 10^{146}$ & \normalsize $2.1\cdot 10^{150}$ & \normalsize $5.1 \cdot 10^{153}$   \\ [1ex]
\hline \\[1ex]

\large $N_{d}$ & \normalsize $42$   & \normalsize $42$  & \normalsize $42$ \\[1ex]
\hline \\[1ex]

\large $N_{t}$ & \normalsize $150$  & \normalsize $150$  & \normalsize $150$  \\[1ex]
\hline\\[1ex]

\end{tabular}
\label{001}
\end{table}

\begin{table}[h]
\caption{$\delta=0.08$}
\centering
\tiny
\begin{tabular}{c llllll}\\[2ex]
\hline\hline \\ [1ex]
\large $c$ & \large $60$ & \large $70$ & \large $80$ & \large $90$ & \large $100$ & \large $110$  \\ [1ex]
\hline\\[1ex]

\large $RMS$ & \normalsize $ 3.8\cdot 10^{-76} $ & \normalsize $2.1\cdot 10^{-78} $ & \normalsize $ 4.8\cdot 10^{-80} $ & \normalsize $ 9.4\cdot 10^{-81} $ & \normalsize $1.2\cdot 10^{-81} $  & \normalsize $1.6\cdot 10^{-82}$   \\ [1ex]
\hline\\[1ex]

\large $COND$ & \normalsize $1.3\cdot 10^{200}$ & \normalsize $2.6\cdot 10^{208}$ & \normalsize $4.0 \cdot 10^{215}$ & \normalsize $8.7 \cdot 10^{221}$ & \normalsize $4.1 \cdot 10^{227}$ & \normalsize $5.5 \cdot 10^{232}$ \\ [1ex]
\hline\\[1ex]

\large $N_{d}$ & \normalsize $63$  & \normalsize $63$  & \normalsize $63$  & \normalsize $63$ & \normalsize $63$  & \normalsize $63$\\[1ex]
\hline\\[1ex]

\large $N_{t}$ & \normalsize $150$  & \normalsize $150$  & \normalsize $150$  & \normalsize $150$ & \normalsize $150$  & \normalsize $150$  \\[1ex]

\hline\hline \\[1ex]

\large $c$ & \large $120$ & \large $130$ & \large $150$ \\ [1ex]
\hline\\[1ex]

\large $RMS$ & \normalsize $ 2.9\cdot 10^{-82} $ & \normalsize $3.4\cdot 10^{-82} $ & \normalsize $ 4.7\cdot 10^{-82} $     \\ [1ex]
\hline \\[1ex]

\large $COND$ & \normalsize $2.7\cdot 10^{237}$ & \normalsize $5.4\cdot 10^{241}$ & \normalsize $2.8\cdot 10^{249}$   \\ [1ex]
\hline \\[1ex]

\large $N_{d}$ & \normalsize $63$   & \normalsize $63$  & \normalsize $63$ \\[1ex]
\hline \\[1ex]

\large $N_{t}$ & \normalsize $150$  & \normalsize $150$  & \normalsize $150$ \\[1ex]
\hline\\[1ex]

\end{tabular}
\label{001}
\end{table}

\clearpage

\begin{table}[t]
\caption{$\delta=0.06$}
\centering
\tiny
\begin{tabular}{c llllll}\\[2ex]
\hline\hline \\ [1ex]
\large $c$ & \large $80$ & \large $90$ & \large $100$ & \large $110$ & \large $120$ & \large $140$  \\ [1ex]
\hline\\[1ex]

\large $RMS$ & \normalsize $ 3.4\cdot 10^{-110} $ & \normalsize $3.3\cdot 10^{-112} $ & \normalsize $ 5.4\cdot 10^{-114} $ & \normalsize $ 1.2\cdot 10^{-115} $ & \normalsize $3.1\cdot 10^{-116} $  & \normalsize $6.2\cdot 10^{-118}$   \\ [1ex]
\hline\\[1ex]

\large $COND$ & \normalsize $6.6\cdot 10^{288}$ & \normalsize $2.0\cdot 10^{297}$ & \normalsize $7.9 \cdot 10^{304}$ & \normalsize $5.9 \cdot 10^{311}$ & \normalsize $1.1 \cdot 10^{318}$ & \normalsize $1.4 \cdot 10^{329}$ \\ [1ex]
\hline\\[1ex]

\large $N_{d}$ & \normalsize $84$  & \normalsize $84$  & \normalsize $84$  & \normalsize $84$ & \normalsize $84$  & \normalsize $84$\\[1ex]
\hline\\[1ex]

\large $N_{t}$ & \normalsize $329$  & \normalsize $329$  & \normalsize $329$  & \normalsize $329$ & \normalsize $329$  & \normalsize $329$  \\[1ex]

\hline\hline \\[1ex]

\large $c$ & \large $160$ & \large $170$ & \large $180$ & \large $200$ \\ [1ex]
\hline\\[1ex]

\large $RMS$ & \normalsize $ 4.9\cdot 10^{-118} $ & \normalsize $4.4\cdot 10^{-118} $ & \normalsize $ 5.3\cdot 10^{-118} $ & \normalsize $ 1.6\cdot 10^{-117} $     \\ [1ex]
\hline \\[1ex]

\large $COND$ & \normalsize $6.0\cdot 10^{338}$ & \normalsize $1.4\cdot 10^{343}$ & \normalsize $1.8 \cdot 10^{347}$ & \normalsize $7.2 \cdot 10^{354}$  \\ [1ex]
\hline \\[1ex]

\large $N_{d}$ & \normalsize $84$   & \normalsize $84$  & \normalsize $84$ & \normalsize $84$  \\[1ex]
\hline \\[1ex]

\large $N_{t}$ & \normalsize $329$  & \normalsize $329$  & \normalsize $329$ & \normalsize $329$  \\[1ex]
\hline\\[1ex]

\end{tabular}
\label{001}
\end{table}

\begin{table}[h]
\caption{$\delta=0.12$}
\centering
\tiny
\begin{tabular}{c llllll}\\[2ex]
\hline\hline \\ [1ex]
\large $c$ & \large $50$ & \large $60$ & \large $70$ & \large $80$ & \large $90$ & \large $100$  \\ [1ex]
\hline\\[1ex]

\large $RMS$ & \normalsize $ 1.7\cdot 10^{-46} $ & \normalsize $3.1\cdot 10^{-47} $ & \normalsize $ 3.9\cdot 10^{-48} $ & \normalsize $ 2.5\cdot 10^{-49} $ & \normalsize $8.6\cdot 10^{-49} $  & \normalsize $5.0\cdot 10^{-49}$   \\ [1ex]
\hline\\[1ex]

\large $COND$ & \normalsize $3.9\cdot 10^{126}$ & \normalsize $1.2\cdot 10^{133}$ & \normalsize $3.6 \cdot 10^{138}$ & \normalsize $2.1 \cdot 10^{143}$ & \normalsize $3.2 \cdot 10^{147}$ & \normalsize $1.6 \cdot 10^{151}$ \\ [1ex]
\hline\\[1ex]

\large $N_{d}$ & \normalsize $42$  & \normalsize $42$  & \normalsize $42$  & \normalsize $42$ & \normalsize $42$  & \normalsize $42$\\[1ex]
\hline\\[1ex]

\large $N_{t}$ & \normalsize $160$  & \normalsize $160$  & \normalsize $160$  & \normalsize $160$ & \normalsize $160$  & \normalsize $160$  \\[1ex]

\hline\hline \\[1ex]

\large $c$ & \large $110$ & \large $120$ & \large $130$ \\ [1ex]
\hline\\[1ex]

\large $RMS$ & \normalsize $ 9.4\cdot 10^{-48} $ & \normalsize $3.9\cdot 10^{-47} $ & \normalsize $ 2.4\cdot 10^{-46} $      \\ [1ex]
\hline \\[1ex]

\large $COND$ & \normalsize $4.0\cdot 10^{154}$ & \normalsize $5.1\cdot 10^{157}$ & \normalsize $3.6 \cdot 10^{160}$   \\ [1ex]
\hline \\[1ex]

\large $N_{d}$ & \normalsize $42$   & \normalsize $42$  & \normalsize $42$ \\[1ex]
\hline \\[1ex]

\large $N_{t}$ & \normalsize $160$  & \normalsize $160$  & \normalsize $160$  \\[1ex]
\hline\\[1ex]

\end{tabular}
\label{001}
\end{table}

\clearpage

\begin{table}[t]
\caption{$\delta=0.08$}
\centering
\tiny
\begin{tabular}{c llllll}\\[2ex]
\hline\hline \\ [1ex]
\large $c$ & \large $50$ & \large $60$ & \large $70$ & \large $80$ & \large $90$ & \large $100$  \\ [1ex]
\hline\\[1ex]

\large $RMS$ & \normalsize $ 3.0\cdot 10^{-72} $ & \normalsize $8.8\cdot 10^{-76} $ & \normalsize $ 4.5\cdot 10^{-78} $ & \normalsize $ 1.3\cdot 10^{-76} $ & \normalsize $2.2\cdot 10^{-80} $  & \normalsize $2.9\cdot 10^{-81}$   \\ [1ex]
\hline\\[1ex]

\large $COND$ & \normalsize $3.2\cdot 10^{191}$ & \normalsize $2.6\cdot 10^{201}$ & \normalsize $5.1 \cdot 10^{209}$ & \normalsize $7.8 \cdot 10^{216}$ & \normalsize $1.7 \cdot 10^{223}$ & \normalsize $8.0 \cdot 10^{228}$ \\ [1ex]
\hline\\[1ex]

\large $N_{d}$ & \normalsize $63$  & \normalsize $63$  & \normalsize $63$  & \normalsize $63$ & \normalsize $63$  & \normalsize $63$\\[1ex]
\hline\\[1ex]

\large $N_{t}$ & \normalsize $320$  & \normalsize $320$  & \normalsize $320$  & \normalsize $320$ & \normalsize $320$  & \normalsize $320$  \\[1ex]

\hline\hline \\[1ex]

\large $c$ & \large $110$ & \large $120$ & \large $130$ & \large $140$  \\ [1ex]
\hline\\[1ex]

\large $RMS$ & \normalsize $ 3.2\cdot 10^{-82} $ & \normalsize $6.8\cdot 10^{-82} $ & \normalsize $ 7.9\cdot 10^{-82} $ & \normalsize $ 1.0\cdot 10^{-81} $  \\ [1ex]
\hline \\[1ex]

\large $COND$ & \normalsize $1.1\cdot 10^{234}$ & \normalsize $5.2\cdot 10^{238}$ & \normalsize $1.1\cdot 10^{243}$ & \normalsize $1.0 \cdot 10^{247}$   \\ [1ex]
\hline \\[1ex]

\large $N_{d}$ & \normalsize $63$   & \normalsize $63$  & \normalsize $63$ & \normalsize $63$  \\[1ex]
\hline \\[1ex]

\large $N_{t}$ & \normalsize $320$  & \normalsize $320$  & \normalsize $320$ & \normalsize $320$  \\[1ex]
\hline\\[1ex]

\end{tabular}
\label{001}
\end{table}

\begin{table}[h]
\caption{$\delta=0.06$}
\centering
\tiny
\begin{tabular}{c llllll}\\[2ex]
\hline\hline \\ [1ex]
\large $c$ & \large $50$ & \large $60$ & \large $70$ & \large $80$ & \large $90$ & \large $100$  \\ [1ex]
\hline\\[1ex]

\large $RMS$ & \normalsize $ 1.1\cdot 10^{-98} $ & \normalsize $3.6\cdot 10^{-103} $ & \normalsize $ 9.1\cdot 10^{-107} $ & \normalsize $ 2.0\cdot 10^{-109} $ & \normalsize $1.9\cdot 10^{-111} $  & \normalsize $3.3\cdot 10^{-113}$   \\ [1ex]
\hline\\[1ex]

\large $COND$ & \normalsize $1.4\cdot 10^{256}$ & \normalsize $1.8\cdot 10^{269}$ & \normalsize $2.3 \cdot 10^{280}$ & \normalsize $9.8 \cdot 10^{289}$ & \normalsize $3.0 \cdot 10^{298}$ & \normalsize $1.2 \cdot 10^{306}$ \\ [1ex]
\hline\\[1ex]

\large $N_{d}$ & \normalsize $84$  & \normalsize $84$  & \normalsize $84$  & \normalsize $84$ & \normalsize $84$  & \normalsize $84$\\[1ex]
\hline\\[1ex]

\large $N_{t}$ & \normalsize $320$  & \normalsize $320$  & \normalsize $320$  & \normalsize $320$ & \normalsize $320$  & \normalsize $320$  \\[1ex]

\hline\hline \\[1ex]

\large $c$ & \large $110$ & \large $120$ & \large $130$ & \large $140$ & \large $150$  \\ [1ex]
\hline\\[1ex]

\large $RMS$ & \normalsize $ 6.7\cdot 10^{-115} $ & \normalsize $1.9\cdot 10^{-115} $ & \normalsize $ 2.0\cdot 10^{-116} $ & \normalsize $ 5.0\cdot 10^{-113} $ & \normalsize $6.0\cdot 10^{-107} $     \\ [1ex]
\hline \\[1ex]

\large $COND$ & \normalsize $8.7\cdot 10^{312}$ & \normalsize $1.6\cdot 10^{319}$ & \normalsize $9.5 \cdot 10^{324}$ & \normalsize $2.1 \cdot 10^{330}$ & \normalsize $2.0 \cdot 10^{335}$  \\ [1ex]
\hline \\[1ex]

\large $N_{d}$ & \normalsize $84$   & \normalsize $84$  & \normalsize $84$ & \normalsize $84$  & \normalsize $84$\\[1ex]
\hline \\[1ex]

\large $N_{t}$ & \normalsize $320$  & \normalsize $320$  & \normalsize $320$ & \normalsize $320$  & \normalsize $320$ \\[1ex]
\hline\\[1ex]

\end{tabular}
\label{001}
\end{table}

\clearpage

\begin{table}[t]
\caption{$\delta=0.03$}
\centering
\tiny
\begin{tabular}{c llllll}\\[2ex]
\hline\hline \\ [1ex]
\large $c$ & \large $40$ & \large $60$ & \large $80$ & \large $100$ & \large $120$ & \large $140$  \\ [1ex]
\hline\\[1ex]

\large $RMS$ & \normalsize $ 5.4\cdot 10^{-189} $ & \normalsize $3.3\cdot 10^{-213} $ & \normalsize $ 1.0\cdot 10^{-229} $ & \normalsize $ 1.5\cdot 10^{-241} $ & \normalsize $2.1\cdot 10^{-250} $  & \normalsize $3.5\cdot 10^{-257}$   \\ [1ex]
\hline\\[1ex]

\large $COND$ & \normalsize $3.4\cdot 10^{479}$ & \normalsize $8.1\cdot 10^{537}$ & \normalsize $2.3 \cdot 10^{579}$ & \normalsize $3.3 \cdot 10^{611}$ & \normalsize $6.3 \cdot 10^{637}$ & \normalsize $1.0 \cdot 10^{660}$ \\ [1ex]
\hline\\[1ex]

\large $N_{d}$ & \normalsize $167$  & \normalsize $167$  & \normalsize $167$  & \normalsize $167$ & \normalsize $167$  & \normalsize $167$\\[1ex]
\hline\\[1ex]

\large $N_{t}$ & \normalsize $700$  & \normalsize $700$  & \normalsize $700$  & \normalsize $700$ & \normalsize $700$  & \normalsize $700$  \\[1ex]

\hline\hline \\[1ex]

\large $c$ & \large $160$ & \large $170$ & \large $180$ & \large $190$ & \large $200$  \\ [1ex]
\hline\\[1ex]

\large $RMS$ & \normalsize $ 1.5\cdot 10^{-262} $ & \normalsize $7.2\cdot 10^{-265} $ & \normalsize $ 2.5\cdot 10^{-269} $ & \normalsize $ 7.2\cdot 10^{-266} $ & \normalsize $5.0\cdot 10^{-259} $     \\ [1ex]
\hline \\[1ex]

\large $COND$ & \normalsize $1.9\cdot 10^{679}$ & \normalsize $1.0\cdot 10^{688}$ & \normalsize $1.8 \cdot 10^{696}$ & \normalsize $1.1 \cdot 10^{704}$ & \normalsize $2.8 \cdot 10^{711}$  \\ [1ex]
\hline \\[1ex]

\large $N_{d}$ & \normalsize $167$   & \normalsize $167$  & \normalsize $167$ & \normalsize $167$  & \normalsize $167$\\[1ex]
\hline \\[1ex]

\large $N_{t}$ & \normalsize $700$  & \normalsize $700$  & \normalsize $700$ & \normalsize $700$  & \normalsize $700$ \\[1ex]
\hline\\[1ex]

\end{tabular}
\label{001}
\end{table}

It is easily seen that in Tables 15-18, the evenly spaced data setting, the optimal values of $c$ go away from the theoretically predicted value 60 as the parameter $\delta$ decreases. It is the same for the scattered data setting, as shown in Tables 19-22. Therefore, one must be careful whenever the bottom of the MN curve tends to be horizontal.
\section{Summary}
We are satisfied with the performance of the MN curve approach to finding the optimal value of the shape parameter, both in the evenly spaced and purely scattered data settings. Although this approach was presented by the author, the foundation built by W.R. Madych and S.A. Nelson plays an important role. Based on this foundation, the author eventually presented a practically useful theory. Hence we name the crucial function MN function, in honor of their outstanding contribution. It is natural to ask whether our theory can be improved. To our regret, the answer probably is `no'. It is already known that algebraic-type error bounds do not reflect the influence of the shape parameter well. As for the exponential-type error bound raised by Madych and Nelson in \cite{MN3}, which applies to scattered data settings, shows the influence of the shape parameter sufficiently only when fill distance is extremely small, making it practically useless. This can be seen in Luh \cite{Lu5}. The improved exponential-type error bound, namely Theorem 1.2 of this paper, shows the influence of the shape parameter sufficiently when fill distance is of reasonable size. In the field of radial basis functions, this kind of exponential-type error bound probably is already optimal, due to the uncertainty principle subject to the condition number, as can be seen in Schaback \cite{Sch}. It means that even if there is an exponential-type error bound which can be used to predict directly the optimal value of the shape parameter and applies to scatered data settings, it may not be better than the approach developed from Theorem 1.2 of this paper.

As for the function space, although $B_{\sigma}$ in Definition 1.1 is quite small, it plays only an intermediate role in the process of the interpolation. We repeatedly pointed out that any function in the Sobolev space, which contains the solutions to many important differential equations, can be interpolated by an $B_{\sigma}$ function with a good error bound, as shown in Narcowich et al. \cite{NW}. Then the $B_{\sigma}$ function can be interpolated by a function in the form of (2) with the same set of data points, also with a good error bound, of which the MN function $MN(c)$ is its essential part. The distance between the Sobolev space function and the RBF interpolator (2) can be handled by triangle inequality. The $B_{\sigma}$ function need not be found explicitly. One only needs to know that it exists. The choice of the parameter $\sigma$ is very flexible. As long as it makes both error bounds small, it is a good choice.


\begin{thebibliography}{99}

\bibitem{MN1}Madych WR. Nelson SA.
{\em Multivariate interpolation and conditionally positive definite function,}
Approx. Theory Appl. 1988;4(4):77-89.


\bibitem{We}Wendland H.
{\em Scattered data approximation,}
Cambridge University Press; 2005.


\bibitem{Fl}Fleming W.
{\em Functions of several variables, second ed.,}
New York: Springer-Verlag; 1977.


\bibitem{Bo}Bos LP.
{\em Bounding the Lebesgue function for Lagrange interpolation in a simplex,}
J. Approx. Theory 1983;38:43-59.


\bibitem{Lu1}Luh L-T.
{\em The equivalence theory of native spaces,}
Approx. Theory Appl. 2001;17(1):76-96.

\bibitem{Lu2}Luh L-T.
{\em The embedding theory of native spaces,}
Approx. Theory Appl. 2001;17(4):90-104.


\bibitem{Lu3}Luh L-T.
{\em On Wu and Schaback's error bound,}
Inter. J. Numer. Methods Appl. 2009;1(2):155-74.


\bibitem{MN2}Madych WR. Nelson SA.
{\em Multivariate interpolation and conditionally positive definite function, II,}
Math. Comp. 1990;54:211-30.

\bibitem{Lu4}Luh L-T.
{\em The mystery of the shape parameter III,}
Appl. Comput. Harmon. Anal. 2016;40:186-199.


\bibitem{MN3}Madych WR. Nelson SA.
{\em Bounds on multivariate polynomials and exponential error estimates for multiquadric interpolation,}
 J. Approx. Theory 1992;70:94-114.


\bibitem{Lu5}Luh L-T.
{\em The mystery of the shape parameter,}
arXiv:1001.5087; 2010.

\bibitem{Sch}Schaback R.
{\em Error estimates and condition numbers for radial basis function interpolation,} 
 Adv. Comput. Math. 1995;3:251-264. 

\bibitem{NW}Narcowich FJ. Ward JD. Wendland H.
{\em  Sobolev error estimates and a Berstein inequality for scattered data interpolation via radial basis functions,}
 Constr. Approx. 2006;20:175-186.



\end{thebibliography}
\end{document}